# Sample Robust Scheduling of Electricity-Gas Systems Under Wind Power Uncertainty

Rong-Peng Liu, *Graduate Student Member*, *IEEE*, Yunhe Hou, *Senior Member*, *IEEE*, Yujia Li, *Student Member*, *IEEE*, Shunbo Lei, *Member*, *IEEE*, Wei Wei, *Senior Member*, *IEEE*, and Xiaozhe Wang, *Senior Member*, *IEEE*

*Abstract*—Bulk integrated electricity and gas systems (IEGSs) introduce complex coupling relations and induce synergistic operation challenges. The growing uncertainty arising from the renewable power generation in the IEGS further aggravates the synergistic problems. Considering the availability of historical wind power generation data, this paper adopts a two-stage sample robust optimization (SRO) model, which is equivalent to the two-stage distributionally robust optimization (DRO) model with a type-∞ Wasserstein ambiguity set, to address the wind power penetrated unit commitment optimal energy flow (UC-OEF) problem for the IEGS. Compared to the equivalent DRO model, the two-stage SRO model can be approximately transformed into a computationally efficient form. Specifically, we employ linear decision rules to simplify the proposed UC-OEF model. Moreover, we further enhance the tractability of the simplified model by exploring its structural features and, accordingly, develop a solution method. Simulation results on two IEGSs validate the effectiveness of the proposed model and solution method.

## Nomenclature

### A. Sets

| | |
|---|---|
| $\mathcal{P}_d/\mathcal{G}_d$ | Set of power/gas loads. |
| $\mathcal{P}_g/\mathcal{G}_g/\mathcal{G}_w$ | Set of coal-fired generators/gas-fired generators/gas wells. |
| $\mathcal{P}_l/\mathcal{G}_l/\mathcal{G}_c$ | Set of power transmission lines/gas passive pipelines/gas compressors (gas active pipelines). |
| $\mathcal{P}_n/\mathcal{G}_n/\mathcal{G}_n^F$ | Set of power nodes/gas nodes/free gas nodes. |
| $\mathcal{P}_w$ | Set of wind farms. |
| $\mathcal{S}$ | Set of scenarios. |
| $\mathcal{T}_T$ | Set of time periods from 1 to $T$. |
| $\mathcal{U}_s^t$ | Uncertainty set. |

### B. Parameters

| | |
|---|---|
| $C_g/C_{w'}$ | Cost of coal-fired generator $g$/gas well $w'$. |
| $CU_g/CD_g$ | Start-up/shut-down cost of generator $g$. |
| $G_n^{min}/G_n^{max}$ | Pressure limit of gas node $n$. |
| $G_{w'}^{min}/G_{w'}^{max}$ | Output limit of gas well $w'$. |
| $ON_g/OF_g$ | Minimum on/off time of generator $g$. |
| $P_d^t/G_d^t$ | Nodal load $d$ in power/gas system at time $t$. |
| $P_g^{min}/P_g^{max}$ | Output limit of generator $g$. |
| $P_w^{min}/P_w^{max}$ | Output limit of wind farm $w$. |
| $P_l/G_l/G_c$ | Transmission limit of power transmission line $l$/gas passive pipeline $l$/gas compressor $c$. |
| $SU_g/SD_g$ | Start-up/shut-down ramp rate of generator $g$. |
| $RU_g/RD_g$ | Ramp-up/ramp-down rate of generator $g$. |
| $W_l$ | Weymouth equation constant of gas passive pipeline $l$. |
| $x_l$ | Reactance of power transmission line $l$. |
| $\alpha_c^{min}/\alpha_c^{max}$ | Ratio limit of gas compressor $c$. |
| $\chi_g$ | Electricity-gas conversion ratio of gas-fired generator $g$. |

### C. Variables

| | |
|---|---|
| $p_g^{t,s}/p_w^{t,s}/g_{w'}^{t,s}$ | Output of generator $g$/wind farm $w$/gas well $w'$ of sample $s$ at time $t$. |
| $p_l^{t,s}/g_l^{t,s}/g_c^{t,s}$ | Power/gas/gas flow through power transmission line $l$/gas passive pipeline $l$/gas compressor $c$ of sample $s$ at time $t$. |
| $x_g^t/u_g^t/v_g^t$ | Binary variable. 1 represents generator $g$ is on/start-up/shut-down at time $t$, and 0 otherwise. |
| $\pi_n^{t,s}$ | Pressure of gas node $n$ of sample $s$ at time $t$. |

This work was supported in part by the National Key R&D Program of China (Technology and Application of Wind Power / Photovoltaic Power Prediction for Promoting Renewable Energy Consumption, 2018YFB0904200) and in part by the Eponymous Complement S&T Program of State Grid Corporation of China (SGLNDKOOKJJS1800266).

R. Liu, Y. Hou, and Y. Li are with the Department of Electrical and Electronic Engineering, the University of Hong Kong, Hong Kong SAR, China, and also with Huazhong University of Science and Technology, Wuhan 430000, China (e-mail: rpliu@eee.hku.hk; yhhou@eee.hku.hk; yjli@eee.hku.hk).

S. Lei is with the School of Science and Engineering, Chinese University of Hong Kong, Shenzhen, Guangdong 518172, China, and also with the Department of Electrical Engineering and Computer Science, University of Michigan, Ann Arbor, MI 48109, USA (email: shunbo.lei@gmail.com).

W. Wei is with the State Key Laboratory of Power Systems, Department of Electrical Engineering, Tsinghua University, Beijing 100084, China (e-mail: wei-wei04@mails.tsinghua.edu.cn).

X. Wang is with the Department of Electrical and Computer Engineering, McGill University, Montreal, QC H3A 0G4, Canada (e-mail: xiaozhe.wang2@mcgill.ca).

## I. Introduction

The mushrooming of the gas-fired power generation spurs the blossom of integrated electricity and gas system (IEGSs), owing to its superiority over the coal-fired generation [1]. However, the interdependency between two systems causes synergistic operation challenges, one of which is the optimal energy flow (OEF) problem [2], [3]. Reference [2] proposes a convexified model to optimize and manage the transmission congestion in the IEGS. Reference [3] proposes a distributed framework to address unit commitment (UC) OEF problems for the multi-area IEGS.

Meanwhile, driven by the eco-friendly requirement and huge profits, renewables are integrated into IEGSs to generate electricity. Consequently, their uncertainty further affects the synergistic operation of IEGSs and leads to *renewable penetrated OEF problems*. Previous research work proposes a variety of modeling paradigms to address these problems. One of the most reliable models which can fully accommodate the renewable fluctuation is the robust optimization (RO) model [4], [5]. In [4], an RO model is adopted to optimize the cost of the wind power penetrated IEGS and ensure its feasibility. Based on an approximated transient matrix-form gas flow model, reference [5] proposes a two-stage robust generation scheduling to cope with the wind uncertainty for IEGSs. As is well known, the RO model is considered to be conservative [6], since it does not incorporate probability information of random variables.

In fact, the probability information of renewable power generation (i.e., random variables) can be estimated by historical data, e.g., historical wind power generation data [7]. Compared with RO models, data-driven modeling methods make the best of data and thus are more applicable to obtaining less conservative operation strategies for renewable penetrated OEF problems. As an approximation of stochastic optimization (SO) models, the sample average approximation (SAA) model is a commonly used data-driven model and has been applied to IEGS optimization problems [8], [9]. Reference [8] employs an SAA model for the day-ahead scheduling of IEGSs under wind uncertainty. Reference [9] adopts an SAA model to optimize the operation cost of wind power penetrated IEGSs and deploy the flexible ramp. Note that the out-of-sample capability of the SAA model may be unpromising when faced with limited data [10]. Even if data is sufficient, the balance between the computational complexity and out-of-sample capability remains a dilemma.

Another data-driven modelling method is the distributionally robust optimization (DRO) model [11]. This model minimizes the expected cost under the worst-case distribution of random variables over an ambiguity set and thus has a better out-of-sample capability than the SAA model and is less conservative than the RO model for power system optimization problems [12]-[14]. Reference [12] adopts a two-stage DRO model for renewable power penetrated optimal energy and reserve dispatch problems. Reference [13] uses the DRO model with an enhanced ambiguity set to address UC problems under wind uncertainty. Reference [14] employs a DRO model to construct new power distribution systems against contingencies. Limited work applies the DRO model to IEGS optimization problems [15]-[17]. Reference [15] proposes a two-stage DRO model with a moment-based ambiguity set to address the day-ahead scheduling of the IEGS considering the demand response and load uncertainty. The same model is adopted in [16] to handle chance-constrained optimal power-gas flow problems under wind uncertainty. According to [18], moment-based ambiguity sets cannot guarantee the convergence of an unknown distribution to its true distribution. In [17], a risk-based DRO model with a Wasserstein-based ambiguity set is proposed to handle wind power penetrated OEF problems. Unfortunately, the Wasserstein ambiguity set-based DRO model is computationally expensive [10].

Recently, reference [10] proposes an equivalent model to the two-stage DRO model with a type-$\infty$ Wasserstein ambiguity set, i.e., the two-stage *sample robust optimization (SRO)* model. Like a combination of SAA and RO models, the two-stage SRO model is constructed by generating multiple uncertainty sets, each of them around one sample. Its first-stage decision aims to accommodate all uncertainties with the lowest mean of the worst-case cost over uncertainty sets. One of the merits of the two-stage SRO model is that its out-of-sample capability depends on both the sample size and the uncertainty budgets in uncertainty sets. Thus, it is capable of handling the situation (with a promising out-of-sample capability) when there are limited historical data available. In addition, the two-stage SRO model can be approximately transformed into a *computationally efficient form* and avoids the computational difficulty of the equivalent DRO model.

Considering the available historical wind power generation data and the out-of-sample capability, we adopt the two-stage SRO model to address the wind power penetrated UC-OEF problem. Furthermore, we enhance its tractability by exploring the structural features of the proposed sample robust UC-OEF model. The main contributions are summarized as follows:

1) This paper adopts a two-stage SRO model to address wind power penetrated UC-OEF problems. Multiple uncertainty sets, each around one wind power generation data (hereinafter called *sample*), are generated in the proposed sample robust UC-OEF model. The objective is to minimize the sum of the (first-stage) UC cost and mean of the (second-stage) worst-case operation cost over these uncertainty sets.

2) This paper simplifies the proposed model using the linear decision rule (LDR). The affine relation between second-stage decision variables and random variables, i.e., wind power outputs, is established via linearized gas transmission equations and free gas nodes. Accordingly, the proposed tri-level UC-OEF model is simplified as a bi-level model (denoted as the LDR-based model hereinafter), and its tractability is enhanced.

3) This paper develops a duality-based solution method to solve the LDR-based model. To enhance the computational efficiency of the original duality-based method, we transform most of the robust constraints into non-robust ones and remove inactive thermal limit constraints by exploring the structural features of the LDR-based model. Testing results demonstrate the effectiveness of the proposed solution method.

In addition to this work, the DRO model and LDR are also adopted in [6], [13], [15], [16], [19], and [20] to address power system-related optimization problems. Differently, references [6] and [13] aim at power system UC problems, and references [19] and [20] address the optimization problems in combined heating and power systems. References [15] and [16] focus on the OEF problem for the IEGS. However, renewable power generation is not considered in [15], and reference [16] ignores the binary UC decisions. Moreover, both of them adopt the moment-based ambiguity set, which cannot guarantee any convergence to the true distribution. One of the most notable differences is that we enhance the tractability of the sample robust UC-OEF model by exploring its structural features, which is never considered in the above references.

In fact, there are other density-based DRO models [21]-[24]. Specifically, reference [21] applies the $L_2$-norm-based DRO model to chance-constrained and two-stage UC problems under renewable uncertainties. Reference [22] proposes $L_1$- and $L_\infty$-norm-based DRO models to address wind power penetrated UC problems. Reference [23] develops a general DRO model with the ambiguity set covering both discrete and continuous distributions by assuming that the relatively complete recourse property holds. Reference [24] proposes a $\phi$-divergence-based chance-constrained DRO model and equivalently reformulates it as a chance-constrained stochastic program. We intend to compare the SRO model with these models in terms of: i) Tractability. The reformulations of all models are tractable, although their tractability varies (mixed-integer second-order cone program (MISOCP) in [21], bi-level (min-max) mixed-integer linear program (MILP) in [22] and [23], general chance-constrained stochastic program in [24], MILP in this work); ii) Parameters. Statistical inference can be used in [21]-[24] to help obtain the parameters in the ambiguity sets, while the uncertainty budget in the SRO is chosen manually; iii) Conservativeness. The conservativeness of all models vanishes as the sample size increases to infinity; iv) Scalability. The size of the reformulations (numbers of decision variables and constraints) grows linearly with the sample size in the SRO model and the number of bins in [22]. As is aforementioned, we can further enhance the tractability of the (SRO-based) MILP reformulation of the wind power penetrated UC-OEF problem by exploring its structural features (see Section III). Thus, the SRO model is adopted in this paper.

The remainder of this paper is organized as follows. Section II presents the two-stage SRO model, sample robust UC-OEF formulation, LDR, and LDR-based model. Section III develops the solution method. Case studies are conducted and reported in Section IV. Section V concludes this paper.

## II. PROBLEM FORMULATION

Firstly, we intend to introduce the scope, assumptions, and simplifications of this work.

1) This paper only considers wind power uncertainty, since compared with wind power forecasting models (e.g., [25]), load forecasting models provide more accurate load profiles [26]. Please refer to [27] for the load uncertainty-based model.

2) This paper considers radial gas networks, while meshed gas networks are not incorporated. According to [17], many gas transmission networks are radial, indicating the practical value of this work.

3) Based on engineering experience, gas flow directions in radial gas networks can be determined in advance and thus are fixed in the proposed model. This assumption is also adopted in [17], [28].

### A. Two-Stage Linear SRO Model

The general two-stage linear SRO model [10] is as follows:

$$\min_{\mathbf{y} \in \mathbb{R}^m} \mathbf{c}_1^T \mathbf{y} + \frac{1}{N} \sum_{i=1}^{N} \max_{\xi_i \in \mathcal{U}_i} \min_{\mathbf{z}_i \in \mathbb{R}^n} \mathbf{c}_2^T \mathbf{z}_i \quad (1a)$$

$$\text{s.t. } \mathbf{A}\mathbf{y} + \mathbf{B}\mathbf{z}_i \leq \mathbf{H}\xi_i + \mathbf{h} \quad \forall \xi_i \in \mathcal{U}_i, \ i \in \{1, \cdots, N\}. \quad (1b)$$

$N$ is the sample size. $\mathbf{y}$ and $\mathbf{z}_i$, $i \in \{1, \cdots, N\}$, are first-stage and second-stage decision variables, respectively. $\xi_i$ are random variables and vary in uncertainty sets $\mathcal{U}_i := \{\xi_i \in \Theta \mid |\xi_i - \bar{\xi}_i| \leq \varepsilon_i\}$, where $\Theta := \{\xi_i \in \mathbb{R}^d \mid \mathbf{Q}\xi_i \leq \mathbf{q}\}$ so that $\mathbb{P}(\xi_i \in \Theta) = 1$ holds for any $i \in \{1, \cdots, N\}$. Namely, the bounds of random variables (not necessarily tight) are known. $\bar{\xi}_i$ are samples of $\xi_i$. $\varepsilon_i$, $\varepsilon_i \geq 0$, are the maximum deviation. $\mathbf{A}$, $\mathbf{B}$, and $\mathbf{H}$ are constant matrices. $\mathbf{c}_1$, $\mathbf{c}_2$, and $\mathbf{h}$ are constant vectors.

The two-stage SRO model (1) reduces to the two-stage SAA model if all $\varepsilon_i = 0$ and reduces to the two-stage RO model if $N = 1$. Since all $\varepsilon_i \geq 0$ in (1), it provides an *upper bound* for the two-stage SAA model and has a promising out-of-sample capability without sacrificing its average performance. More importantly, it is *asymptotically optimal*, indicating that the gap between the optimums of (1) and the two-stage SAA model tends to zero when $N \to \infty$ and $\varepsilon_i \to 0$ [10]. Note that model (1) becomes a (general) two-stage mixed-integer linear SRO model if its first-stage decision variable $\mathbf{y}$ consists of integers or mixed integers.

### B. Sample Robust UC-OEF Model

The sample robust UC-OEF model is as follows:

$$\min_{\{\mathbf{x},\mathbf{u},\mathbf{v}\}} \sum_{t \in \mathcal{T}_T} \sum_{g \in \mathcal{P}_g} C_g^0 x_g^t + CU_g u_g^t + CD_g v_g^t$$

$$+ \frac{1}{|\mathcal{S}|} \sum_{s \in \mathcal{S}} \sum_{t \in \mathcal{T}_T} \max_{\mathbf{p}_W^{t,s} \in \mathcal{U}_s^t} \min_{\mathbf{z}^{t,s}} \sum_{g \in \mathcal{P}_g} C_g p_g^{t,s} + \sum_{w' \in \mathcal{G}_w} C_{w'} g_{w'}^{t,s} \quad (2)$$

s.t. $x_g^t - x_g^{t-1} = u_g^t - v_g^t \quad \forall g \in \mathcal{P}_g \cup \mathcal{G}_g,\ t \in \mathcal{T}_T$ (3a)

$$ON_g u_g^t \leq \sum_{k=t}^{t+ON_g-1} x_g^k \quad \forall g \in \mathcal{P}_g \cup \mathcal{G}_g,\ t \in \mathcal{T}_{T-ON_g+1} \quad (3b)$$

$$OF_g v_g^t \leq \sum_{k=t}^{t+OF_g-1} (1-x_g^k) \quad \forall g \in \mathcal{P}_g \cup \mathcal{G}_g,\ t \in \mathcal{T}_{T-OF_g+1} \quad (3c)$$

$$\sum_{k=t}^{T} u_g^k \leq \sum_{k=t}^{T} x_g^k \quad \forall g \in \mathcal{P}_g \cup \mathcal{G}_g,\ t \in \mathcal{T}_T \setminus \mathcal{T}_{T-ON_g} \quad (3d)$$

$$\sum_{k=t}^{T} v_g^k \leq \sum_{k=t}^{T} (1-x_g^k) \quad \forall g \in \mathcal{P}_g \cup \mathcal{G}_g,\ t \in \mathcal{T}_T \setminus \mathcal{T}_{T-OF_g} \quad (3e)$$

$$P_g^{\min} x_g^t \leq p_g^{t,s} \leq P_g^{\max} x_g^t \quad \forall g \in \mathcal{P}_g \cup \mathcal{G}_g,\ t \in \mathcal{T}_T,\ s \in \mathcal{S} \quad (4a)$$

$$-RD_g x_g^t - SD_g v_g^t \leq p_g^{t,s} - p_g^{t-1,s} \leq RU_g x_g^{t-1} + SU_g u_g^t$$
$$\forall g \in \mathcal{P}_g \cup \mathcal{G}_g,\ t \in \mathcal{T}_T,\ s \in \mathcal{S} \quad (4b)$$

$$-P_l \leq \sum_{k} \beta_{lk} \Big( \sum_{g \in \mathcal{P}_g \cup \mathcal{G}_g} A_{kg}^G p_g^{t,s} + \sum_{w \in \mathcal{P}_w} A_{kw}^W p_w^{t,s} - \sum_{d \in \mathcal{P}_d} A_{kd}^D P_d^t \Big) \leq P_l$$
$$\forall \mathbf{p}_W^{t,s} \in \mathcal{U}_s^t,\ l \in \mathcal{P}_l,\ t \in \mathcal{T}_T,\ s \in \mathcal{S} \quad (4c)$$

$$\sum_{g \in \mathcal{P}_g \cup \mathcal{G}_g} p_g^{t,s} + \sum_{w \in \mathcal{P}_w} p_w^{t,s} = \sum_{d \in \mathcal{P}_d} P_d^t \quad \forall \mathbf{p}_W^{t,s} \in \mathcal{U}_s^t,\ t \in \mathcal{T}_T,\ s \in \mathcal{S} \quad (4d)$$

$$G_{w'}^{\min} \leq g_{w'}^{t,s} \leq G_{w'}^{\max} \quad \forall w' \in \mathcal{G}_w,\ t \in \mathcal{T}_T,\ s \in \mathcal{S} \quad (5a)$$

$$G_n^{\min} \leq \pi_n^{t,s} \leq G_n^{\max} \quad \forall n \in \mathcal{G}_n,\ t \in \mathcal{T}_T,\ s \in \mathcal{S} \quad (5b)$$

$$g_l^{t,s} = W_l \sqrt{((\pi_{m(l)}^{t,s})^2 - (\pi_{n(l)}^{t,s})^2)} \quad \forall l \in \mathcal{G}_l,\ t \in \mathcal{T}_T,\ s \in \mathcal{S} \quad (5c)$$

$$\alpha_c^{\min} \pi_{n(c)}^{t,s} \leq \pi_{m(c)}^{t,s} \leq \alpha_c^{\max} \pi_{n(c)}^{t,s} \quad \forall c \in \mathcal{G}_c,\ t \in \mathcal{T}_T,\ s \in \mathcal{S} \quad (5d)$$

$$0 \leq g_l^{t,s} \leq G_l \quad \forall l \in \mathcal{G}_l,\ t \in \mathcal{T}_T,\ s \in \mathcal{S} \quad (5e)$$

$$0 \leq g_c^{t,s} \leq G_c \quad \forall c \in \mathcal{G}_c,\ t \in \mathcal{T}_T,\ s \in \mathcal{S} \quad (5f)$$

$$\sum_{w' \in \mathcal{G}_w} g_{w'(n)}^{t,s} + \sum_{l_1 \in \mathcal{G}_l} g_{l_1(n)}^{t,s} - \sum_{l_2 \in \mathcal{G}_l} g_{l_2(n)}^{t,s} + \sum_{c_1 \in \mathcal{G}_c} g_{c_1(n)}^{t,s} - \sum_{c_2 \in \mathcal{G}_c} g_{c_2(n)}^{t,s}$$
$$= \sum_{d \in \mathcal{G}_d} G_{d(n)}^t + \sum_{g \in \mathcal{G}_g} \chi_g p_{g(n)}^{t,s} \quad \forall n \in \mathcal{G}_n,\ t \in \mathcal{T}_T,\ s \in \mathcal{S}. \quad (5g)$$

Objective (2) aims to minimize the sum of the first-stage UC cost (non-load, start-up, and shut-down cost) and the second-stage re-dispatch cost (output cost from coal-fired generators and gas wells). $C_g^0$ is the non-load cost of generator $g$. For any $t \in \mathcal{T}_T$ and $s \in \mathcal{S}$, $\mathcal{U}_s^t := \{\mathbf{p}_W^{t,s} \in \Xi \mid \sum_{w \in \mathcal{P}_w} (|p_w^{t,s} - \bar{p}_w^{t,s}|) \leq$ $\varepsilon_s^t \sum_{w \in \mathcal{P}_w} \bar{p}_w^{t,s}$ }, and $\Xi := \{ \mathbf{p}_W^{t,s} \in \mathbb{R}^{|\mathcal{P}_w|} \mid P_w^{\min} \leq p_w^{t,s} \leq P_w^{\max}, w \in \mathcal{P}_w \}$, where $\bar{p}_w^{t,s}$ are wind power generation samples. Vector $\mathbf{p}_W^{t,s} = col(p_w^{t,s})$, $w \in \mathcal{P}_w$, $\mathbf{p}_W^s = col(\mathbf{p}_W^{t,s})$, $t \in \mathcal{T}_T$ and $\mathbf{p}_W = col(\mathbf{p}_W^s)$, $s \in \mathcal{S}$. We can obtain $\mathbf{p}_G^s, \mathbf{g}_{W'}^s, \boldsymbol{\pi}_N^s, \mathbf{g}_L^s, \mathbf{g}_C^s, \mathbf{p}_G, \mathbf{g}_{W'}, \boldsymbol{\pi}_N, \mathbf{g}_L$, and $\mathbf{g}_C$ in the similar way. Hereinafter, all (new) vector variables are generated using the same rule. Vectors $\mathbf{x}$, $\mathbf{u}$, and $\mathbf{v}$ consist of first-stage UC decision variables $x_g^t, u_g^t$, and $v_g^t$, respectively, and vector $\mathbf{z}^{t,s} = col(\mathbf{p}_G^{t,s}, \mathbf{g}_{W'}^{t,s}, \boldsymbol{\pi}_N^{t,s}, \mathbf{g}_L^{t,s}, \mathbf{g}_C^{t,s})$, i.e., second-stage decision variables.

Power system constraints are denoted by (3)-(4). Constraint (3a) enforces the on and off states of both coal-fired and gas-fired generators, and constraints (3b)-(3e) regulate their minimum on and off time. Constraint (4a) restricts the output range of a generator, and its ramping capability is limited by constraint (4b). Constraint (4c), which is based on the direct current (DC) power flow model, states the power transmission line capacity. $A_{kg}^G, A_{kw}^G$, and $A_{kd}^G$ are the elements in incidence matrices for generators, wind farms, and power loads, respectively. $\beta_{lk}$ is the power transfer distribution factor [29]. Equation (4d) is the power balance constraint.

Gas system constraints are denoted by (5). Constraint (5a) enforces the output capacity of a gas well. The range of nodal pressure is bounded by (5b). Equation (5c), known as the Weymouth equation, is widely leveraged to model the gas flow transported in long-distance gas passive pipelines [2], [8], where $m(l)$ and $n(l)$ are the nodes connected by pipeline $l$. Constraint (5d) is a simplified gas compressor model [15] connecting gas nodes $m(c)$ and $n(c)$. Constraints (5e) and (5f) restrict the gas flow in a gas passive pipeline and a gas compressor, respectively. Equation (5g) balances the supply and demand at each gas node, where $w'(n), l_1(n), l_2(n), c_1(n), c_2(n), d(n)$, and $g(n)$ represent the gas well, gas inflows and outflows of the gas passive pipeline and compressor, gas loads, and gas-fired generator connected by node $n$, respectively.

The two-stage sample robust UC-OEF model (2)-(5) is computationally intractable due to the nonconvex Weymouth equation (5c). In this paper, the Weymouth equation is linearized by the Taylor series approximation [2], [30], i.e.,

$$g_l^{t,s} = W_l \sqrt{((\pi_{m(l)}^{t,s})^2 - (\pi_{n(l)}^{t,s})^2)} \approx K_{m(l)}^{t,s} \pi_{m(l)}^{t,s} - K_{n(l)}^{t,s} \pi_{n(l)}^{t,s}$$
$$\forall l \in \mathcal{G}_l, t \in \mathcal{T}_T, s \in \mathcal{S}, \quad (6)$$

where $K_{m(l)}^{t,s}, K_{m(l)}^{t,s} = W_l \bar{\pi}_{m(l)}^{t,s} / \sqrt{\bar{\pi}_{m(l)}^{t,s} - \bar{\pi}_{n(l)}^{t,s}}$, and $K_{n(l)}^{t,s}, K_{n(l)}^{t,s} = W_l \bar{\pi}_{n(l)}^{t,s} / \sqrt{\bar{\pi}_{m(l)}^{t,s} - \bar{\pi}_{n(l)}^{t,s}}$, are linearization points, which are critical for the accuracy of the approximation. $\bar{\pi}_{m(l)}^{t,s}$ and $\bar{\pi}_{n(l)}^{t,s}$ are constants. This work adopts the method proposed in [17] to obtain linearization points, and the procedures are summarized as follows: i) set $\varepsilon_s^t = 0$ for all $t \in \mathcal{T}_T$ and $s \in \mathcal{S}$ in model (2)-(5); ii) replace the Weymouth equation (5c) with its second-order cone form; iii) solve the deterministic model specified in Steps i) and ii) by the sequential second-order cone program-based Algorithm 1 in [17]; iv) set $\bar{\pi}_{m(l)}^{t,s}$ and $\bar{\pi}_{n(l)}^{t,s}$ to corresponding gas pressure values (specified in Step iii) and obtain all linearization points in (6).

### C. Free Gas Node

For radial gas networks, we introduce the concept of the *free gas node*. At a free gas node, the optimal solution for (2)-(5) is not unique. It is the compressor model (5d) that results in the non-uniqueness. An Illustrative example is presented in Fig. 1. The radial gas network in the left part of Fig. 1 consists of $i$ compressors, $k$-$i$ gas passive pipelines, and $k$+1 nodes. By removing compressors, this network is separated into $i$+1 gas subnetworks (shown in the right part of Fig. 1). Without loss of generality, we assume that the reference gas node belongs to $GN_{k+1}$ (hereinafter called reference subnetwork). Nodes 1 to $i$ are free gas nodes. This example indicates that: i) a free gas node is connected by at least one compressor; ii) except for the reference subnetwork, each subnetwork has one and only one free gas node; iii) there is no free gas node in the reference subnetwork. Namely, the number of free gas nodes equals to the number of gas compressors. Based on the above analysis, we propose Method 1 to locate free gas nodes systematically.

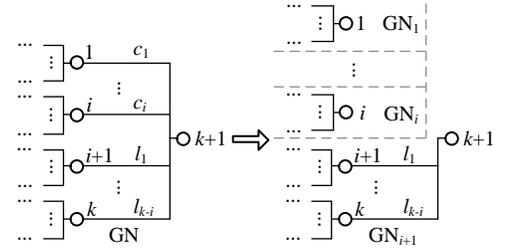

Fig. 1. Illustrative diagram of free gas nodes in a radial gas network.

**Method 1**: Free gas node locating
1: Initialize $\mathcal{G}_n^F \leftarrow \emptyset$, $\mathcal{G}_c^F \leftarrow \mathcal{G}_c$
2: remove all compressors from the radial gas network
3: **do**
4:   new network $\leftarrow$ merging the reference subnetwork with one of its neighboring subnetworks by adding $c_i$, $c_i \in \mathcal{G}_c$
5:   $\mathcal{G}_n^F \leftarrow \mathcal{G}_n^F \cup \{n(c)\}$, $n(c) \notin$ reference subnetwork
6:   $\mathcal{G}_c^F \leftarrow \mathcal{G}_c^F \setminus \{c_i\}$
7:   reference subnetwork $\leftarrow$ new network
8: **until** $\mathcal{G}_c^F = \emptyset$
9: output: $\mathcal{G}_n^F$

The purpose of introducing free gas nodes is to denote other gas variables by gas well outputs and free gas node pressures (please refer to Appendix for details). This is beneficial to reduce the number of variables and, more importantly, build rational linear decision rules (see the next subsection).

### D. Linear Decision Rule Based Model

It is still challenging to solve the tri-level sample robust UC-OEF model (2)-(5b), (5d)-(6). We consider simplifying this model by the LDR proposed in [27]. For any $t \in \mathcal{T}_T$ and $s \in \mathcal{S}$, the output of each generator, the output of each gas well, and the pressure of each free gas node are regarded as affine functions of the sum of wind farm outputs, respectively, i.e.,

$$p_g^{t,s}(\mathbf{p}_W^{t,s}) = r_g^{t,s} + R_g^{t,s} \sum_{w \in \mathcal{P}_w} p_w^{t,s} \quad \forall g \in \mathcal{P}_g \cup \mathcal{G}_g \quad (7a)$$

$$g_{w'}^{t,s}(\mathbf{p}_W^{t,s}) = s_{w'}^{t,s} + S_{w'}^{t,s} \sum_{w \in \mathcal{P}_w} p_w^{t,s} \quad \forall w' \in \mathcal{G}_w \quad (7b)$$

$$\pi_n^{t,s}(\mathbf{p}_W^{t,s}) = o_n^{t,s} + O_n^{t,s} \sum_{w \in \mathcal{P}_w} p_w^{t,s} \quad \forall n \in \mathcal{G}_n^F, \quad (7c)$$

where $r_g^{t,s}$, $R_g^{t,s}$, $s_{w'}^{t,s}$, $S_{w'}^{t,s}$, $o_n^{t,s}$, and $O_n^{t,s}$ are (coefficient) variables of the LDR (7). Note that affine relations (7a) and (7b) hold (due to (4d) and (5g)), provided a fixed (first-stage) UC decision, and (7c) is only an approximation. Using Taylor series model (6), free gas node, and LDR (7), constraints (4) and (5) are reformulated as

$$P_g^{\min} x_g^t \leq r_g^{t,s} + R_g^{t,s} \sum_{w \in \mathcal{P}_w} p_w^{t,s} \leq P_g^{\max} x_g^t$$
$$\forall \mathbf{p}_W^{t,s} \in \mathcal{U}_s^t, \, g \in \mathcal{P}_g \cup \mathcal{G}_g, \, t \in \mathcal{T}_T, \, s \in \mathcal{S} \quad (8a)$$

$$-RD_g x_g^t - SD_g v_g^t \leq r_g^{t,s} + R_g^{t,s} \sum_{w \in \mathcal{P}_w} p_w^{t,s} - (r_g^{t-1,s} + R_g^{t-1,s} \sum_{w \in \mathcal{P}_w} p_w^{t-1,s})$$
$$\leq RU_g x_g^{t-1} + SU_g u_g^t$$
$$\forall \mathbf{p}_W^{t,s} \in \mathcal{U}_s^t, \, \mathbf{p}_W^{t-1,s} \in \mathcal{U}_s^{t-1}, \, g \in \mathcal{P}_g \cup \mathcal{G}_g, \, t \in \mathcal{T}_T, \, s \in \mathcal{S} \quad (8b)$$

$$-P_l \leq \sum_k \beta_{lk} \big( \sum_{g \in \mathcal{P}_g \cup \mathcal{G}_g} A_{kg}^G (r_g^{t,s} + R_g^{t,s} \sum_{w \in \mathcal{P}_w} p_w^{t,s}) + \sum_{w \in \mathcal{P}_w} A_{kw}^W p_w^{t,s}$$
$$- \sum_{d \in \mathcal{P}_d} A_{kd}^D P_d^t \big) \leq P_l \quad \forall \mathbf{p}_W^{t,s} \in \mathcal{U}_s^t, \, l \in \mathcal{P}_l, \, t \in \mathcal{T}_T, \, s \in \mathcal{S} \quad (8c)$$

$$\sum_{g \in \mathcal{P}_g \cup \mathcal{G}_g} r_g^{t,s} = \sum_{d \in \mathcal{P}_d} P_d^{t,s}, \quad \sum_{g \in \mathcal{P}_g \cup \mathcal{G}_g} R_g^{t,s} = -1 \quad \forall t \in \mathcal{T}_T, \, s \in \mathcal{S} \quad (8d)$$

$$G_{w'}^{\min} \leq s_{w'}^{t,s} + S_{w'}^{t,s} \sum_{w \in \mathcal{P}_w} p_w^{t,s} \leq G_{w'}^{\max}$$
$$\forall \mathbf{p}_W^{t,s} \in \mathcal{U}_s^t, \, w' \in \mathcal{G}_w, \, t \in \mathcal{T}_T, \, s \in \mathcal{S} \quad (9a)$$

$$G_n^{\min} \leq o_n^{t,s} + O_n^{t,s} \sum_{w \in \mathcal{P}_w} p_w^{t,s} \leq G_n^{\max}$$
$$\forall \mathbf{p}_W^{t,s} \in \mathcal{U}_s^t, \, n \in \mathcal{G}_n^F, \, t \in \mathcal{T}_T, \, s \in \mathcal{S}$$
$$G_n^{\min} \leq q_n(\mathbf{b}^{t,s}) \leq G_n^{\max} \quad \forall \mathbf{p}_W^{t,s} \in \mathcal{U}_s^t, \, n \in \mathcal{G}_n \setminus \mathcal{G}_n^F, \, t \in \mathcal{T}_T, \, s \in \mathcal{S} \quad (9b)$$

$$\begin{cases} \alpha_c^{\min}(o_{n(c)}^{t,s} + O_{n(c)}^{t,s} \sum_{w \in \mathcal{P}_w} p_w^{t,s}) \leq q_{m(c)}(\mathbf{b}^{t,s}) \\ \qquad \leq \alpha_c^{\max}(o_{n(c)}^{t,s} + O_{n(c)}^{t,s} \sum_{w \in \mathcal{P}_w} p_w^{t,s}) \quad \text{if } n \in \mathcal{G}_n^F \\ \alpha_c^{\min} q_{n(c)}(\mathbf{b}^{t,s}) \leq o_{m(c)}^{t,s} + O_{m(c)}^{t,s} \sum_{w \in \mathcal{P}_w} p_w^{t,s} \leq \alpha_c^{\max} q_{n(c)}(\mathbf{b}^{t,s}) \\ \qquad \qquad \qquad \qquad \qquad \qquad \qquad \text{if } m \in \mathcal{G}_n^F \end{cases}$$
$$\forall \mathbf{p}_W^{t,s} \in \mathcal{U}_s^t, \, c \in \mathcal{G}_c, \, t \in \mathcal{T}_T, \, s \in \mathcal{S} \quad (9c)$$

$$0 \leq \sum_k Q_{lk}^P \big( \sum_{w' \in \mathcal{G}_w} B_{kw'}^{W'} (s_{w'}^{t,s} + S_{w'}^{t,s} \sum_{w \in \mathcal{P}_w} p_w^{t,s}) - \sum_{d \in \mathcal{G}_d} B_{kd}^{D'} G_d^t$$
$$- \sum_{g \in \mathcal{G}_g} B_{kg}^{G'} \chi_g (r_g^{t,s} + R_g^{t,s} \sum_{w \in \mathcal{P}_w} p_w^{t,s})) \leq G_l$$
$$\forall \mathbf{p}_W^{t,s} \in \mathcal{U}_s^t, \, l \in \mathcal{G}_l, \, t \in \mathcal{T}_T, \, s \in \mathcal{S} \quad (9d)$$

$$0 \leq \sum_k Q_{ck}^P \big( \sum_{w' \in \mathcal{G}_w} B_{kw'}^{W'} (s_{w'}^{t,s} + S_{w'}^{t,s} \sum_{w \in \mathcal{P}_w} p_w^{t,s}) - \sum_{d \in \mathcal{G}_d} B_{kd}^{D'} G_d^t$$
$$- \sum_{g \in \mathcal{G}_g} B_{kg}^{G'} \chi_g (r_g^{t,s} + R_g^{t,s} \sum_{w \in \mathcal{P}_w} p_w^{t,s})) \leq G_c$$
$$\forall \mathbf{p}_W^{t,s} \in \mathcal{U}_s^t, \, c \in \mathcal{G}_c, \, t \in \mathcal{T}_T, \, s \in \mathcal{S} \quad (9e)$$

$$\sum_{w' \in \mathcal{G}_w} s_{w'}^{t,s} = \sum_{g \in \mathcal{G}_g} \chi_g r_g^{t,s} + \sum_{d \in \mathcal{G}_d} G_d^t, \quad \sum_{w' \in \mathcal{G}_w} S_{w'}^{t,s} = \sum_{g \in \mathcal{G}_g} \chi_g R_g^{t,s}$$
$$\forall t \in \mathcal{T}_T, \, s \in \mathcal{S}, \quad (9f)$$

where $\mathbf{b}^{t,s} = \{\mathbf{r}_G^{t,s}, \mathbf{R}_G^{t,s}, \mathbf{s}_{W'}^{t,s}, \mathbf{S}_{W'}^{t,s}, \mathbf{o}_N^{t,s}, \mathbf{O}_N^{t,s}, \mathbf{p}_W^{t,s}\}$ and

$$q_n(\mathbf{b}^{t,s}) = \sum_k Q_{nk}^{NP} \big( \sum_{w' \in \mathcal{G}_w} B_{kw'}^{W'} (s_{w'}^{t,s} + S_{w'}^{t,s} \sum_{w \in \mathcal{P}_w} p_w^{t,s}) - \sum_{d \in \mathcal{G}_d} B_{kd}^{D'} G_d^t$$
$$- \sum_{g \in \mathcal{G}_g} B_{kg}^{G'} \chi_g (r_g^{t,s} + R_g^{t,s} \sum_{w \in \mathcal{P}_w} p_w^{t,s})) - \sum_{k \in \mathcal{G}_n^F} Q_{nk}^{NF} (o_k^{t,s} + O_k^{t,s} \sum_{w \in \mathcal{P}_w} p_w^{t,s}).$$

Please refer to Appendix for constants $Q_{lk}^P$, $Q_{ck}^P$, $B_{kw'}^{W'}$, $B_{kd}^{D'}$, $B_{kg}^{G'}$, $Q_{nk}^{NP}$, and $Q_{nk}^{NF}$. Accordingly, model (2)-(5) is reformulated as

$$\min_{\substack{\{\mathbf{x}, \mathbf{u}, \mathbf{v}, \mathbf{r}_G, \mathbf{R}_G, \\ \mathbf{s}_{W'}, \mathbf{S}_{W'}, \mathbf{o}_N, \mathbf{O}_N\}}} \sum_{t \in \mathcal{T}_T} \sum_{g \in \mathcal{P}_g} C_g^0 x_g^t + CU_g u_g^t + CD_g v_g^t + \frac{1}{|\mathcal{S}|} \sum_{s \in \mathcal{S}} \sum_{t \in \mathcal{T}_T} \max_{\mathbf{p}_W^{t,s} \in \mathcal{U}_s^t}$$
$$\big( (\sum_{g \in \mathcal{P}_g} C_g (r_g^{t,s} + R_g^{t,s} \sum_{w \in \mathcal{P}_w} p_w^{t,s})) + \sum_{w' \in \mathcal{G}_w} C_{w'} (s_{w'}^{t,s} + S_{w'}^{t,s} \sum_{w \in \mathcal{P}_w} p_w^{t,s}) \big) \quad (10a)$$

s.t. constraints (3), (8) and (9). (10b)

By means of the LDR, the proposed tri-level sample robust UC-OEF model (2)-(5b), (5d)-(6) is simplified as the two-level LDR-based model (10). Generally, the LDR-based model is a *conservative estimation* of its original model [6], [27]. Thus, model (10) may induce a higher operation cost and suboptimal (first-stage) UC decision. Differently, reference [10] proposes a promising asymptotically optimal property for the general multi-policy LDR-based two-stage SRO model (see Theorem 2 in [10]). By applying it to (10), we conclude that under some mild conditions, the optimal objective value and optimal UC decision of the LDR-based model (10) converge to those of the sample robust UC-OEF model (2)-(5b), (5d)-(6) when $N \to \infty$ and all $\varepsilon_i \to 0$. This property indicates that theoretically, we can decrease the conservativeness of (10) (caused by the LDR) to any range by increasing the sample size. Namely, model (10) can achieve a flexible balance between the computational cost and quality of the UC decision according to actual needs.

### III. Solution Method

The model (10) can be directly reformulated as a single-level optimization problem using the duality theory. However, this method introduces mass new (dual) variables and constraints which induce computational intractability, even for a moderate-sized system [27]. This section simplifies model (10) before the duality-based reformulation by exploring the structural features of polyhedral uncertainty sets $\mathcal{U}_s^t$ and the LDR.

#### A. Robust Constraint Transformation

Define the constraints containing random variables as *robust constraints* and the constraints without any random variables as *non-robust constraints*. In model (10), constraints (8a)-(8c) and (9a)-(9e) are robust constraints, and constraints (3), (8d), and (9f) are non-robust constraints. We rewrite the model (10) as the following compact form:

$$\min_{\mathbf{z}'} \mathbf{c}_1^T \mathbf{z}'' + \frac{1}{|\mathcal{S}|} \sum_{s \in \mathcal{S}} \max_{\mathbf{p}_W^s \in \mathcal{U}_s} \mathbf{c}_2^T (\mathbf{y}^s + \mathbf{Y}^s \mathbf{p}_W^s) \quad (11a)$$

s.t. $\mathbf{A}_0 \mathbf{z}' - \mathbf{a}_0 \geq 0$ (11b)

$\mathbf{A} \mathbf{z}'' + \mathbf{B}(\mathbf{y}^s + \mathbf{Y}^s \mathbf{p}_W^s) - (\mathbf{H} \mathbf{p}_W^s + \mathbf{h}_0) \geq 0 \quad \forall \mathbf{p}_W^s \in \mathcal{U}_s, \, s \in \mathcal{S}.$ (11c)

$\mathbf{z}' = \{\mathbf{x}, \mathbf{u}, \mathbf{v}, \mathbf{y}^1, \cdots, \mathbf{y}^{|\mathcal{S}|}, \mathbf{Y}^1, \cdots, \mathbf{Y}^{|\mathcal{S}|}\}$. $\mathbf{z}'' = col(\mathbf{x}, \mathbf{u}, \mathbf{v})$. $\mathbf{y}^s = col(\mathbf{r}_G^s, \mathbf{s}_{W'}^s, \mathbf{o}_N^s)$ and $\mathbf{Y}^s = \{\mathbf{R}_G^s, \mathbf{S}_{W'}^s, \mathbf{O}_N^s\}$. $\mathbf{A}_0$, $\mathbf{A}$, $\mathbf{B}$, and $\mathbf{H}$ are constant matrices. $\mathbf{c}_1$, $\mathbf{c}_2$, $\mathbf{a}_0$, and $\mathbf{h}_0$ are constant vectors. $\mathcal{U}_s = \prod_{\otimes t \in \mathcal{T}_T} \mathcal{U}_s^t$. $\prod_\otimes$ denotes Cartesian product. Constraints (11b) and (11c) refer to non-robust constraints ((3), (8d), and (9f)) and robust constraints ((8a)-(8c) and (9a)-(9e)), respectively.

For model (11), it is proved that the optimal value of $\mathbf{p}_W^s$ (denoted as $\mathbf{p}_W^{s*}$) belongs to the set $ext(\mathcal{U}_s)$ [27], where $ext(\mathcal{U}_s)$ consists of all extreme points in $\mathcal{U}_s$. One of the methods to solve

(11) is the constraint generation [27], i.e., iteratively adding extreme points to a set for each robust constraint in (11c) until finding a solution such that constraints (11c) hold for any $\mathbf{p}_W^s \in \mathcal{U}_s$ and $s \in \mathcal{S}$. We use this idea to simplify (11).

Reference [27] points out that: i) the robust constraint (8a) is equivalent to replacing the uncertainty set $\mathcal{U}_s^t$ with the set $\{\mathbf{p}_{W,\max}^{t,s\,*}, \mathbf{p}_{W,\min}^{t,s\,*}\}$ (since the optimal value for $\sum_{w\in\mathcal{P}_w} p_w^{t,s}$ in (8a) must be one of the extreme points of the polyhedral uncertainty set $\mathcal{U}_s^t$), where $\mathbf{p}_{W,\max}^{t,s\,*}$ and $\mathbf{p}_{W,\min}^{t,s\,*}$ are specified by (12); ii) the robust constraint (8b) is equivalent to replacing uncertainty sets $\mathcal{U}_s^t$ and $\mathcal{U}_s^{t-1}$ with sets $\{\mathbf{p}_{W,\max}^{t,s\,*}, \mathbf{p}_{W,\min}^{t,s\,*}\}$ and $\{\mathbf{p}_{W,\max}^{t-1,s\,*}, \mathbf{p}_{W,\min}^{t-1,s\,*}\}$, respectively. Thus, constraints (8a) and (8b) are transformed into non-robust constraints. We extend this idea to the gas system model. Specifically, the robust constraints (9a)-(9e) are equivalently transformed into non-robust constraints by replacing uncertainty sets $\mathcal{U}_s^t$ with sets $\{\mathbf{p}_{W,\max}^{t,s\,*}, \mathbf{p}_{W,\min}^{t,s\,*}\}$, respectively. This transformation greatly reduces the number of (newly generated) dual variables when reformulating (11) using the duality theory, since only the robust constraints need to be handled by the duality-based method. By far, thermal limit constraint (8c) remains the only kind of robust constraints.

$$\mathbf{p}_{W,\max}^{t,s\,*} = \arg\max_{\mathbf{p}_W^{t,s}} \sum_{w\in\mathcal{P}_w} p_w^{t,s} \quad \text{s.t.} \quad \mathbf{p}_W^{t,s} \in \mathcal{U}_s^t \quad \forall t \in \mathcal{T}_T, s \in \mathcal{S} \quad (12a)$$

$$\mathbf{p}_{W,\min}^{t,s\,*} = \arg\min_{\mathbf{p}_W^{t,s}} \sum_{w\in\mathcal{P}_w} p_w^{t,s} \quad \text{s.t.} \quad \mathbf{p}_W^{t,s} \in \mathcal{U}_s^t \quad \forall t \in \mathcal{T}_T, s \in \mathcal{S}. \quad (12b)$$

### B. Inactive Thermal Limit Elimination

In practice, thermal limit constraints for some power transmission lines are inactive [6]. Accordingly, we develop Method 2 to eliminate inactive thermal limit constraints.

---
**Method 2**: Inactive thermal limit elimination
1: Initialize: $\mathcal{P}_{l,\max}^R \leftarrow \emptyset$, $\mathcal{P}_{l,\min}^R \leftarrow \emptyset$
2: **for** $l \in \mathcal{G}_l$ **do**
3: $\quad \mathbf{p}_{L,\max}^* \leftarrow$ optimal $p_l^{t,s}$ of (13); $\mathbf{p}_{L,\min}^* \leftarrow$ optimal $p_l^{t,s}$ of (14)
4: $\quad \mathcal{P}_{l,\max}^R \leftarrow \mathcal{P}_{l,\max}^R \cup \{(l,t,s)\}$ for any $p_{l,\max}^{t,s\,*} > \mathrm{P}_l$; $\mathcal{P}_{l,\min}^R \leftarrow \mathcal{P}_{l,\min}^R \cup \{(l,t,s)\}$ for any $p_{l,\min}^{t,s\,*} < -\mathrm{P}_l$
5: **end**
6: output: $\mathcal{P}_{l,\max}^R$ and $\mathcal{P}_{l,\min}^R$

---

For any $l \in \mathcal{P}_l$, models (13) and (14) are as follows:

$$\arg\max_{\{\mathbf{p}_g, \mathbf{p}_l, \mathbf{p}_w\}} \sum_{s\in\mathcal{S}} \sum_{t\in\mathcal{T}_T} p_l^{t,s} \quad (13a)$$

s.t. $0 \le p_g^{t,s} \le \mathrm{P}_g^{\max} \quad \forall g \in \mathcal{P}_g \cup \mathcal{G}_g, t \in \mathcal{T}_T, s \in \mathcal{S} \quad (13b)$

$$\sum_{g\in\mathcal{P}_g\cup\mathcal{G}_g} p_g^{t,s} + \sum_{w\in\mathcal{P}_w} p_w^{t,s} = \sum_{d\in\mathcal{P}_d} \mathrm{P}_d^t, \; \mathbf{p}_W^{t,s} \in \mathcal{U}_s^t \quad \forall t \in \mathcal{T}_T, s \in \mathcal{S} \quad (13c)$$

$$p_l^{t,s} = \sum_k \beta_{lk} \left( \sum_{g\in\mathcal{P}_g\cup\mathcal{G}_g} \mathrm{A}_{kg}^G p_g^{t,s} + \sum_{w\in\mathcal{P}_w} \mathrm{A}_{kw}^W p_w^{t,s} - \sum_{d\in\mathcal{P}_d} \mathrm{A}_{kd}^D \mathrm{P}_d^t \right)$$
$$\forall t \in \mathcal{T}_T, s \in \mathcal{S}, \quad (13d)$$

$$\arg\min_{\{\mathbf{p}_g, \mathbf{p}_l, \mathbf{p}_w\}} \sum_{s\in\mathcal{S}} \sum_{t\in\mathcal{T}_T} p_l^{t,s} \quad (14a)$$

s.t. constraints (13b)-(13d). $\quad (14b)$

The basic idea of Method 2 is that the optimal values $\mathbf{p}_{L,\max}^*$ and $\mathbf{p}_{L,\min}^*$ specified by (13) and (14) are the relaxation of their optimal values specified by (11) (referred to as $\mathbf{p}_L^*$). Namely, the relation $\mathbf{p}_{L,\min}^* \le \mathbf{p}_L^* \le \mathbf{p}_{L,\max}^*$ always holds. Therefore, Method 2 can effectively eliminate inactive thermal limit constraints and reduce the number of robust constraints in (11).

### C. Duality-Based Reformulation

After the above steps, most of the robust constraints in (11) are transformed into non-robust constraints or eliminated, and the equivalent form of (11) becomes

$$\min_{\mathbf{z}''} \mathbf{c}_1^T \mathbf{z}'' + \frac{1}{|\mathcal{S}|} \sum_{s\in\mathcal{S}} \sum_{t\in\mathcal{T}_T} \max_{\mathbf{p}_W^{t,s} \in \mathcal{U}_s^t} \mathbf{c}_2'^T (\mathbf{y}^{t,s} + \mathbf{Y}^{t,s} \mathbf{p}_W^{t,s}) \quad (15a)$$

s.t. $\mathbf{A}_0' \mathbf{z}' - \mathbf{a}_0' \ge 0 \quad (15b)$

$$\mathbf{e}_l^T (\mathbf{B}'(\mathbf{y}^{t,s} + \mathbf{Y}^{t,s} \mathbf{p}_W^{t,s}) - (\mathbf{H}' \mathbf{p}_W^{t,s} + \mathbf{h}_0')) \ge 0$$
$$\forall \mathbf{p}_W^{t,s} \in \mathcal{U}_s^t, \{l,t,s\} \in \mathcal{P}_{l,\max}^R \cup \mathcal{P}_{l,\min}^R, \quad (15c)$$

where $\mathbf{y}^{t,s} = col(\mathbf{r}_G^{t,s}, \mathbf{s}_{W'}^{t,s}, \mathbf{o}_N^{t,s})$ and $\mathbf{Y}^{t,s} = \{\mathbf{R}_G^{t,s}, \mathbf{S}_{W'}^{t,s}, \mathbf{O}_N^{t,s}\}$. $\mathbf{A}_0'$, $\mathbf{B}'$, and $\mathbf{H}'$ are constant matrices. $\mathbf{c}_2'$, $\mathbf{a}_0'$, and $\mathbf{h}_0'$ are constant vectors. The duality-based reformulation of (15) is

$$\min_{\mathbf{z}' \cup \{\lambda, \gamma, \theta, \eta\}} \mathbf{c}_1^T \mathbf{z}'' + \frac{1}{|\mathcal{S}|} \sum_{s\in\mathcal{S}} \sum_{t\in\mathcal{T}_T} (\mathbf{c}_2'^T (\mathbf{y}^{t,s} + \mathbf{Y}^{t,s} \overline{\mathbf{p}}_W^{t,s}) + \varepsilon_s^t \lambda^{t,s}$$
$$+ (\mathbf{U} \overline{\mathbf{p}}_W^{t,s} + \mathbf{u}_0)^T \gamma^{t,s}) \quad (16a)$$

s.t. $\mathbf{A}_0' \mathbf{z}' - \mathbf{a}_0' \ge 0 \quad (16b)$

$$\|(\mathbf{Y}^{t,s})^T \mathbf{c}_2' + \mathbf{U}^T \gamma^{t,s}\|_\infty \le \lambda^{t,s} \quad \forall t \in \mathcal{T}_T, s \in \mathcal{S} \quad (16c)$$

$$\mathbf{e}_l^T (\mathbf{B}'(\mathbf{y}^{t,s} + \mathbf{Y}^{t,s} \overline{\mathbf{p}}_W^{t,s}) - (\mathbf{H}' \overline{\mathbf{p}}_W^{t,s} + \mathbf{h}_0')) \ge \varepsilon_s^t \theta_l^{t,s}$$
$$+ (\mathbf{U} \overline{\mathbf{p}}_W^{t,s} + \mathbf{u}_0)^T \eta_l^{t,s} \quad \forall \{l,t,s\} \in \mathcal{P}_{l,\max}^R \cup \mathcal{P}_{l,\min}^R \quad (16d)$$

$$\|(\mathbf{H}' - \mathbf{B}' \mathbf{Y}^{t,s})^T \mathbf{e}_l + \mathbf{U}^T \eta_l^{t,s}\|_\infty \le \theta_l^{t,s}$$
$$\forall \{l,t,s\} \in \mathcal{P}_{l,\max}^R \cup \mathcal{P}_{l,\min}^R \quad (16e)$$

$$\lambda^{t,s} \ge 0, \; \gamma^{t,s} \ge 0 \quad \forall t \in \mathcal{T}_T, s \in \mathcal{S} \quad (16f)$$

$$\theta_l^{t,s} \ge 0, \; \eta_l^{t,s} \ge 0 \quad \forall \{l,t,s\} \in \mathcal{P}_{l,\max}^R \cup \mathcal{P}_{l,\min}^R, \quad (16g)$$

where $\lambda^{t,s}$, $\gamma^{t,s}$, $\theta_l^{t,s}$, and $\eta_l^{t,s}$ are dual variables. $\|\cdot\|_\infty$ denotes $L_\infty$-norm. $\mathbf{U}\mathbf{p}_W^{t,s} + \mathbf{u}_0$ is the compact representation of set $\Xi$. $\mathbf{U}$ is the constant matrix, and $\mathbf{u}_0$ is a constant vector. $\mathbf{1}$ is all one vector. $\mathbf{e}_l$ is the $l$-th column of the identity matrix. The model (16) is a single-level MILP and can be solved by commercial solvers. For ease of reading, the framework of the proposed solution method is shown in Fig. 2.

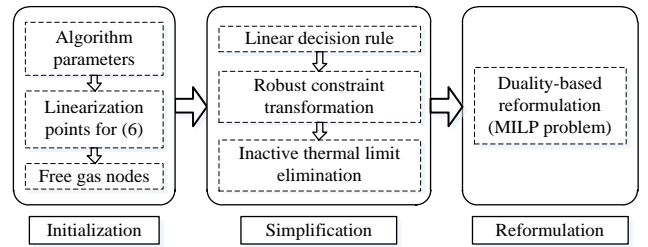

Fig. 2 Framework of the proposed solution method.

## IV. CASE STUDY

In this section, we test the proposed model and solution method on two IEGSs, i.e., the integrated 6-node electricity and 7-node gas system (IEGS-6-7) and the integrated 118-node electricity and 20-node gas system (IEGS-118-20). Please refer to [31] for the system data. Wind power generation samples are partly taken from the Western Wind Data Set of National

Renewable Energy Laboratory (NREL) [7] and partly generated randomly based on the real data. Both tests are coded by MATLAB R2018b with YALMIP toolbox and performed on a PC with an i5-6500 @3.2GHz CPU and a 16 GB memory. The MILP and MISOCP are solved by Gurobi 8.1.1, and the mixed-integer nonlinear program (MINLP) is solved by SCIP 5.0.1.

### A. Integrated 6-Node Electricity and 7-Node Gas System

The IEGS-6-7 has two gas-fired generators. Two wind farms are added to power node 4, and their total installed capacity is 122MW, which accounts for 24.8% of the installed capacity of the IEGS. $T$ is set to 24 with a one-hour interval between two consecutive periods. Convergence thresholds for MILP, MISOCP, and MINLP are set to 1e-4.

*1) Accuracy of the Taylor Series Model*. This part tests the accuracy of the Taylor series model (6), and the results are shown in Table I. M1 (short for Model 1) and M2 (short for Model 2) consist of (2)-(5) and (2)-(5b), (5d)-(6), respectively. For this test, all $\varepsilon_s^t$, $t \in \mathcal{T}_T$, in M1 and M2 are set to 0. $|\mathcal{S}|$ is set to 1, as the sample size has little influence on the test but greatly affects the convergence of M1. Under these settings, M1 and M2 reduce to MINLP and MILP, respectively. PL, GL, and WP represent power loads, gas loads, and wind power generation, respectively, and the budget is their increased percentage. For example, PL under +5% budget means increasing 5% power loads at each power node for $t \in \mathcal{T}_T$. Error = (|Cost of M1−Cost of M2|/Cost of M1)*100%. INF is short for "infeasible". The linearization points in M2 are specified under 0% budget and keep unchanged all the time.

TABLE I
ACCURACY OF THE TAYLOR SERIES MODEL

|  | Budget | +1% | +5% | +10% | +15% | +20% |
|---|---|---|---|---|---|---|
| WP | Cost of M1 | 578170.9 | 574593.7 | 570557.3 | 567642.3 | 564326.8 |
|  | Cost of M2 | 578170.9 | 574595.0 | 570559.3 | 567650.3 | 564332.7 |
|  | Error | <0.001% | <0.001% | <0.001% | 0.001% | 0.001% |
| PL | Cost of M1 | 584732.4 | 607314.7 | 640854.3 | INF | INF |
|  | Cost of M2 | 584738.7 | 607347.7 | 640900.5 | INF | INF |
|  | Error | 0.001% | 0.005% | 0.007% | - | - |
| GL | Cost of M1 | 587595.5 | 604018.1 | 622935.1 | 645226.9 | 670010.4 |
|  | Cost of M2 | 587599.5 | 604045.4 | 623004.4 | 645356.9 | 670194.0 |
|  | Error | <0.001% | 0.004% | 0.011% | 0.020% | 0.027% |

According to this table, M2 performs well under different WP budgets, and all errors are smaller than or equal to 0.001%. Considering the limited wind farm capacity, we also test M2 under different PL and GL budgets. The maximum error is smaller than 0.03%, which is still satisfactory. In addition, M2 correctly reports infeasibility under +15% and +20% PL budgets. In fact, existing short-term load forecasting techniques provide comparatively accurate day-ahead load profiles (no more than 5%) [26], and M2 can easily address this uncertainty budget (errors ≤ 0.005%). Overall, the Taylor series model achieves promising accuracy in approximating the Weymouth equation, especially under the wind power uncertainty.

*2) LDR-Based Model*. This part tests the LDR-based model (10) (denoted as M3). Generally, the LDR-based model is a conservative estimation of its original model [6], [27], and the optimum of M3 (denoted as $opt_{M3}^*$) is larger than or equal to the optimum of M2 (denoted as $opt_{M2}^*$), where M2 consists of (2)-(5b), (5d)-(6). In addition, $opt_{M2}^*$ is an upper bound for the optimum of $\widetilde{M2}$ (denoted as $\widetilde{opt}_{M2}^*$), where $\widetilde{M2}$ consists of (2)-(5b), (5d)-(6) with all $\varepsilon_s^t = 0$ for $t \in \mathcal{T}_T$ and $s \in \mathcal{S}$. Namely, we can choose any legal $\varepsilon_s^t$ for M2 while all $\varepsilon_s^t$ in $\widetilde{M2}$ must equal to 0. We have

$$0 \leq opt_{M3}^* - opt_{M2}^* \leq opt_{M3}^* - \widetilde{opt}_{M2}^*.$$

Theoretically, $opt_{M3}^* - \widetilde{opt}_{M2}^*$ converges to zero when all $\varepsilon_s^t$ in M3 tend to zero. Hence, we can evaluate the LDR-based model by observing the convergence of $opt_{M3}^* - \widetilde{opt}_{M2}^*$ as $\varepsilon_s^t \to 0$. Tests are conducted under different numbers of scenarios ($|\mathcal{S}| = $ 1, 10, 50, and 100), and the results are shown in Fig. 3. Red and blue lines denote $\widetilde{opt}_{M2}^*$ and $opt_{M3}^*$, respectively. Budgets decrease from ten to zero. In Fig. 3, all the $opt_{M3}^*$ in (a)-(d) converge to corresponding $\widetilde{opt}_{M2}^*$ almost linearly when $\varepsilon_s^t \to 0$. Nonlinear fluctuations are caused by different UC decisions. This example shows that the LDR-based model (M3) can be used to approximate the original model (M2).

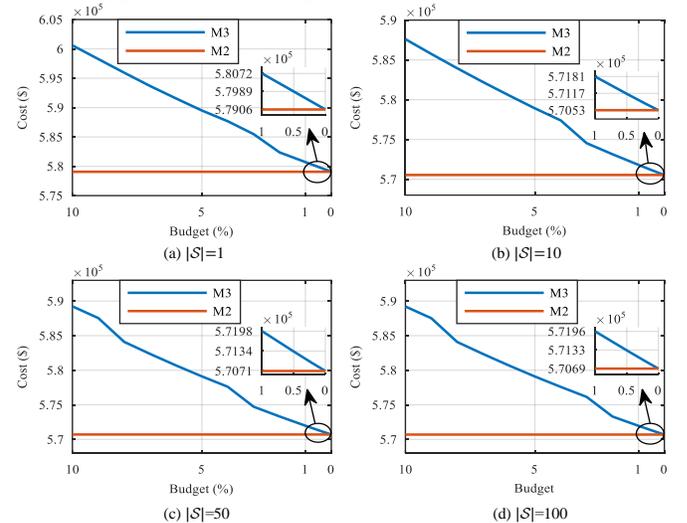

Fig. 3 LDR-based model under different wind power uncertainty budgets.

*3) Computation Time*. This part investigates the effectiveness of the proposed solution method, and the results are listed in Table II. CG and DBM denote constraint generation and duality-based methods, respectively. RCT and TLE refer to robust constraint transformation (Section III.A) and inactive thermal limit elimination (Section III.B) methods, respectively. "+" is the combination of these methods. The computation time denotes the CPU time spending on solving (10) under different uncertainty budgets $\varepsilon_s^t$. The number of scenarios $|\mathcal{S}|$ is set to 100. The CG (mentioned in Section III.A) cannot obtain a solution within one hour. The DBM outperforms the CG and obtains feasible solutions. In the sequel, we only focus on the DBM. The RCT reduces more than half of the computation time of the DBM, since most of the robust constraints are transformed into non-robust constraints and do not need to be reformulated as the dual form. The DBM+RCT+TLE performs best. In order to further test the TLE, thermal limits are relaxed to 125% and 150% of their original values, respectively. Results are shown in the last two rows of Table II. Computation time is further reduced as thermal limit constraints becoming looser. Note that the DBM is the solution method proposed in [10] to solve the two-stage type-∞ Wasserstein-based DRO model and is adopted in this paper to solve the sample robust UC-OEF model. Namely, the computation time of the DBM and DBM+RCT+TLE denotes the computational performance of the two-stage type-∞ Wasserstein-based DRO model in [10] and the two-stage SRO model in this paper for the wind power penetrated UC-OEF problem. Testing results demonstrate the effectiveness of the proposed solution method (DBM+RCT+TLE) in enhancing the computational efficiency of the solution method in [10] (DBM).

TABLE II
COMPUTATION TIME (S)

| Method | $\varepsilon_s^t$=1% | $\varepsilon_s^t$=5% | $\varepsilon_s^t$=10% | $\varepsilon_s^t$=15% | $\varepsilon_s^t$=20% |
|---|---|---|---|---|---|
| CG | >3600 | >3600 | >3600 | >3600 | >3600 |
| DBM | 685 | 703 | 749 | 706 | 731 |
| DBM+RCT | 233 | 237 | 211 | 266 | 237 |
| DBM+RCT+TLE | 140 | 174 | 142 | 160 | 168 |
| DBM+RCT+TLE (125%) | 66 | 79 | 84 | 125 | 143 |
| DBM+RCT+TLE (150%) | 40 | 43 | 57 | 81 | 99 |

*4) Comparison with Other Models.* In this part, the out-of-sample capability of the SRO model is tested and compared with the SAA model. Table III shows the results. One hundred wind power generation samples are newly generated and used to test the out-of-sample performance of SAA and SRO models. "Out-of-sample" in this table refers to the percentage of the newly generated samples which can be accommodated by the first-stage UC decision of the corresponding model. Thermal limits are relaxed to 150% to relieve the computational burden. The SAA model obtained by 100 samples, i.e., $|\mathcal{S}|$=100, performs poorly and can only address 53% of the new samples. By extending samples to 1000, the out-of-sample performance of the (new) SAA model with $|\mathcal{S}|$=1000 is improved. As a comparison, the SRO model provides a more powerful and flexible out-of-sample capability (by adjusting the uncertainty budget according to the actual demand) with less computation time and a smaller sample size. Namely, in order to have a comparable out-of-sample capability to the SRO model, the SAA model requires more samples (than the SRO model), which dramatically increases the computational cost. These results indicate that the SRO model is more applicable to wind power penetrated UC-OEF problems.

TABLE III
COMPARISON BETWEEN SRO AND SAA MODELS

| Model | SAA | | SRO ($|\mathcal{S}|$=100) | |
|---|---|---|---|---|
| | $|\mathcal{S}|$=100 | $|\mathcal{S}|$=1000 | $\varepsilon_s^t$=1% | $\varepsilon_s^t$=5% |
| Out-of-sample | 53.0% | 74.0% | 74.0% | 82.0% |
| Computation time (s) | 4 | 122 | 40 | 43 |

Furthermore, the SRO model is compared with the RO model. Table IV shows the results. The mean of all samples ($|\mathcal{S}|$=100) is used as the center of the uncertainty set in the RO model, and its uncertainty budget is the sum of i) and ii), where i) is the maximum absolute value of the differences between samples and the uncertainty center, and ii) is the budget in the SRO model. The LDR (7) and duality-based method are adopted to reformulate the RO model into a tractable MILP. Cost refers to the objective value (10). From Table IV, we know that the SRO model provides less conservative and feasible solutions, while the RO model is infeasible under $\varepsilon_s^t$=10% due to its conservativeness. In short, the SRO model reduces the conservativeness of the RO model effectively for the UC-OEF problem. We notice that the computation time of the SRO model nearly doubles that of the RO model for both $\varepsilon_s^t$=1% and $\varepsilon_s^t$=5%. We do not obtain the computation time for the RO model under $\varepsilon_s^t$=10% due to its infeasibility. In the next subsection, we will validate the scalability of the proposed SRO model and solution method using a larger IEGS.

TABLE IV
COMPARISON BETWEEN SRO AND RO MODELS

| Model | $\varepsilon_s^t$=1% | | $\varepsilon_s^t$=5% | | $\varepsilon_s^t$=10% | |
|---|---|---|---|---|---|---|
| | RO | SRO | RO | SRO | RO | SRO |
| Cost (*10^5 $) | 5.351 | 5.249 | 6.313 | 6.016 | INF | 6.086 |
| Computation time (s) | 19 | 40 | 20 | 43 | - | 42 |

*B. Integrated 118-Node Electricity and 20-Node Gas System*

The IEGS-118-20 has fifteen gas-fired generators. Five wind farms are added to the power system [31]. The total installed capacity of wind farms accounts for 12.3% of the installed capacity of the IEGS. $T$ is set to 24. Convergence thresholds for MILP and MISOCP are set to 1e-3.

*1) Computation Time.* This part tests the proposed solution method on the IEGS-118-20. Please refer to Section IV.A 3) for details about CG, DBM, and DBM+RCT+TLE. Table V shows the results. The number of scenarios $|\mathcal{S}|$ is set to 50. Both the CG and DBM fail to obtain feasible solutions within three hours. By contrast, the proposed solution method (DBM+RCT+TLE) successfully solves the proposed model in less than an hour under all uncertainty budgets, verifying its tractability and scalability. Besides, we observe that it usually takes more time to solve the SRO model with a larger uncertainty budget $\varepsilon_s^t$ (using the same solution method). This phenomenon also exists in the IEGS-6-7 (see Table II). The plausible reason is that the search space is enlarged by a larger $\varepsilon_s^t$, resulting in a longer computation time.

TABLE V
COMPUTATION TIME (S)

| Method | $\varepsilon_s^t$=5% | $\varepsilon_s^t$=10% | $\varepsilon_s^t$=15% | $\varepsilon_s^t$=20% | $\varepsilon_s^t$=30% |
|---|---|---|---|---|---|
| CG | >10800 | >10800 | >10800 | >10800 | >10800 |
| DBM | >10800 | >10800 | >10800 | >10800 | >10800 |
| DBM+RCT+TLE | 1756 | 1573 | 2442 | 2756 | 2830 |

*2) Comparison with Other Models.* This part compares the proposed model with SAA and RO models, and test results are shown in Table VI. One thousand newly generated samples are used to test the out-of-sample performance. $|\mathcal{S}|$ is set to 50, and $\varepsilon_s^t$ of the RO model are 20%. We adopt the same method (as shown in Section IV.A 4)) to construct the uncertainty set for the RO model and obtain its MILP reformulation. According to test results, the proposed SRO model has a better out-of-sample capability than the SAA model and is effective in reducing the conservativeness of the RO model for the UC-OEF problem. Moreover, we obtain comparable computation time of the RO model (1964 seconds) and the SRO model (shown in the last row of Table V). In fact, the RO model can be regarded as a special SRO model, in which there is only one scenario with a large uncertainty budget. From this perspective, the computation time of the SRO model should be much longer than that of the RO model. It is the proposed method, i.e., the robust constraint transformation and the inactive thermal limit elimination, that significantly reduces the number of the new constraints and variables generated by the duality-based method. This test fully validates the scalability of the proposed solution method and shows its applicable potential to large-scale IEGSs.

TABLE VI
COMPARISON BETWEEN SRO, SAA, AND RO MODELS

| Model | SAA | RO | SRO | | |
|---|---|---|---|---|---|
| | | | $\varepsilon_s^t$=5% | $\varepsilon_s^t$=10% | $\varepsilon_s^t$=20% |
| Cost (*10$^6$ \$) | 7.313 | 7.941 | 7.364 | 7.415 | 7.527 |
| Out-of-sample | 68.1% | 100% | 90.6% | 92.3% | 98.7% |

## V. CONCLUSION

This paper adopts a two-stage SRO model to address the UC-OEF problem under wind power uncertainty. We employ the LDR to simplify the proposed sample robust UC-OEF model and enhance its tractability. Specifically, we develop robust constraint transformation and inactive thermal limit elimination methods to reduce the computational burden of the duality-based solution method. Test results demonstrate that: i) the LDR-based model performs well in approximating the original model; ii) the proposed duality-based solution method is effective in improving the computational efficiency; iii) the proposed model has a better out-of-sample capability than the SAA model and is less conservative than the RO model for wind power penetrated UC-OEF problems.

## APPENDIX

In the Appendix, we consider a single-scenario single-period model, i.e., $T$=1 and $|\mathcal{S}|$=1, in which the radial gas network has $|\mathcal{G}_n|$ nodes. For the sake of simplicity, all superscripts regarding $t$ and $s$ are removed.

Similar to the power injection, define gas injection as

$$g_{\text{inj},n} = \sum_{w' \in \mathcal{G}_w} g_{w'(n)} - \sum_{d \in \mathcal{G}_d} g_{d(n)} - \sum_{g \in \mathcal{G}_g} \chi_g p_{g(n)}$$

$$\forall n \in \mathcal{G}_n. \quad (17)$$

Its vector form is

$$\mathbf{g}_{\text{inj}} = \mathbf{B}^W \mathbf{g}_{W'} - \mathbf{B}^D \mathbf{G}_D - \mathbf{B}^G \mathbf{B}^\chi \mathbf{g}_G, \quad (18)$$

where $\mathbf{g}_{W'}$=$col(g_{w'})$, $w' \in \mathcal{G}_w$, $\mathbf{G}_D$=$col(G_d)$, $d \in \mathcal{G}_d$, and $\mathbf{g}_G$= $col(p_g)$, $g \in \mathcal{G}_g$. $\mathbf{B}^W$, $\mathbf{B}^D$, and $\mathbf{B}^G \mathbf{B}^\chi$ are incidence matrices for gas wells, gas loads, and gas-fired generators, respectively. The vector form of gas nodal balance equation (5g) is

$$\mathbf{g}_{\text{inj}} = \mathbf{B}^P [\mathbf{g}_L^T \ \mathbf{g}_C^T]^T, \quad (19)$$

where $\mathbf{g}_L$=$col(g_l)$, $l \in \mathcal{G}_l$ and $\mathbf{g}_C$= $col(g_c)$, $c \in \mathcal{G}_c$. $\mathbf{B}^P$ is the incidence matrix. In a radial gas network with $|\mathcal{G}_n|$ nodes, the number of gas pipelines (including both passive pipelines and compressors) is $|\mathcal{G}_n|$ −1. Without loss of generality, we assume that node 1 is the reference node. Let $\mathbf{g}'_{\text{inj}}$=$(g_{\text{inj},2},\cdots,g_{\text{inj},|\mathcal{G}_n|})^T$ and

$$\mathbf{g}'_{\text{inj}} = \mathbf{B}^{W'} \mathbf{g}_W - \mathbf{B}^{D'} \mathbf{G}_D - \mathbf{B}^{G'} \mathbf{B}^\chi \mathbf{g}_G, \quad (20)$$

where $\mathbf{B}^{W'}$, $\mathbf{B}^{D'}$, and $\mathbf{B}^{G'}$ are obtained by deleting first rows of $\mathbf{B}^W$, $\mathbf{B}^D$, and $\mathbf{B}^G$, respectively. Equation (19) is equivalent to

$$\mathbf{g}'_{\text{inj}} = \mathbf{B}^{P'} [\mathbf{g}_L^T \ \mathbf{g}_C^T]^T, \quad (21a)$$

$$\sum_{w \in \mathcal{G}_w} g_w = \sum_{d \in \mathcal{G}_d} G_d^t + \sum_{g \in \mathcal{G}_g} \chi_g p_g. \quad (21b)$$

Matrix $\mathbf{B}^{P'}$ is obtained by deleting the first row of $\mathbf{B}^P$ and is a $(|\mathcal{G}_n| -1) \times (|\mathcal{G}_n| -1)$ matrix. It is straightforward to prove that $\mathbf{B}^{P'}$ is invertible. Thus, $\mathbf{g}_L$ and $\mathbf{g}_C$ are denoted by $\mathbf{g}_W$ and $\mathbf{g}_G$, i.e.,

$$[\mathbf{g}_L^T \ \mathbf{g}_C^T]^T = \mathbf{Q}^P \mathbf{g}'_{\text{inj}}, \quad (22)$$

where $\mathbf{Q}^P$=$(\mathbf{B}^{P'})^{-1}$.

In addition, the vector form of the linearized Weymouth equation (6) is equivalent to

$$\mathbf{g}_L = \mathbf{B}^F \boldsymbol{\pi}_{N_1} + \mathbf{B}^N \boldsymbol{\pi}_{N_2}, \quad (23)$$

where $\boldsymbol{\pi}_{N_1}$=$col(\pi_n)$, $n \in \mathcal{G}_n^F \cup \{1\}$, and $\boldsymbol{\pi}_{N_2}$= $col(\pi_n)$, $n \in \mathcal{G}_n \backslash (\mathcal{G}_n^F \cup \{1\})$. $\mathbf{B}^F$ and $\mathbf{B}^N$ are incidence matrices. It is straightforward to prove that $\mathbf{B}^F$ is invertible. Thus, $\boldsymbol{\pi}_{N_2}$ are denoted by $\mathbf{g}_W$, $\mathbf{g}_G$, and $\boldsymbol{\pi}_{N_1}$, i.e.,

$$\boldsymbol{\pi}_{N_2} = \mathbf{Q}^N (\mathbf{Q}_L^P \mathbf{g}'_{\text{inj}} - \mathbf{B}^F \boldsymbol{\pi}_{N_1}), \quad (24)$$

where $\mathbf{Q}^N$=$(\mathbf{B}^N)^{-1}$. $\mathbf{Q}_L^P$ consists of the first $|\mathcal{G}_l|$ rows of $\mathbf{Q}^P$. Let $\mathbf{Q}^{NP}$=$\mathbf{Q}^N \mathbf{Q}_L^P$ and $\mathbf{Q}^{NF}$=$\mathbf{Q}^N \mathbf{B}^F$. Equation (24) is simplified as

$$\boldsymbol{\pi}_{N_2} = \mathbf{Q}^{NP} \mathbf{g}'_{\text{inj}} - \mathbf{Q}^{NF} \boldsymbol{\pi}_{N_1}. \quad (25)$$